\documentclass{amsart}

\usepackage{amsmath,amssymb,amsthm,latexsym}

\newtheorem{theorem}{Theorem}
\newcommand{\bt}{\begin{theorem}}
\newcommand{\et}{\end{theorem}}
\newtheorem{lemma}{Lemma}
\newcommand{\bl}{\begin{lemma}}
\newcommand{\el}{\end{lemma}}
\newtheorem{corollary}{Corollary}
\newcommand{\bc}{\begin{corollary}}
\newcommand{\ec}{\end{corollary}}
\newcommand{\bconj}{\begin{conjecture}}
\newcommand{\econj}{\end{conjecture}}
\newtheorem{problem}{Problem}
\newcommand{\bprob}{\begin{problem}}
\newcommand{\eprob}{\end{problem}}
\newcommand{\beq}{\begin{equation}}
\newcommand{\eeq}{\end{equation}}
\newcommand{\benum}{\begin{enumerate}}
\newcommand{\eenum}{\end{enumerate}}
\newcommand{\N}{\ensuremath{ \mathbf N }}
\newcommand{\Z}{\ensuremath{\mathbf Z}}

\newcommand{\mca}{\ensuremath{ \mathcal A}}
\newcommand{\mcb}{\ensuremath{ \mathcal B}}

\newcommand{\mcr}{\ensuremath{ \mathcal R}}

\newcommand{\mcx}{\ensuremath{ \mathcal X}}

\newcommand{\mbx}{\ensuremath{ \mathbf x}}
\newcommand{\mby}{\ensuremath{ \mathbf y}}

\DeclareMathOperator{\card}{\text{card}}

\newcommand{\bmat}{\left(\begin{matrix}}
\newcommand{\emat}{\end{matrix}\right)}
\newcommand{\bsmallmat}{\left(\begin{smallmatrix}}
\newcommand{\esmallmat}{\end{smallmatrix}\right)}
\DeclareMathOperator{\qand}{\quad\text{and}\quad}

\DeclareMathOperator{\Tet}{\text{Tet}}

\title[Sumset sizes in additive number theory]
{Triangular and tetrahedral number differences of sumset sizes 
in additive number theory}
\author{Melvyn B.  Nathanson}
\address{Department of Mathematics\\Lehman College (CUNY)\\
Bronx, NY 10468}
\email{melvyn.nathanson@lehman.cuny.edu}

\date{\today}


\thanks{Supported in part by  PSC-CUNY Research Award Program grant 66197-00 54.}

\begin{document}

\begin{abstract}
The study of sums of finite sets of integers has mostly concentrated on sets 
with very small sumsets (Freiman's theorem and related work) 
and on sets with very large sumsets (Sidon sets and $B_h$-sets).  
This paper considers the full range of sumset sizes of finite sets 
of integers and an unexpected 
pattern (related to the triangular and tetrahedral numbers) 
that appears in the distribution of 
popular sumset sizes of sets of size 4. 
\end{abstract}

\maketitle

\section{The sumset size problem}

The \emph{integer interval} $[u,v]$ is the set of all integers $n$ 
such that $u \leq n \leq v$.

The $h$-fold sum of a set $A$ of integers is the set $hA$ of all sums of $h$ 
not necessarily distinct elements of $A$: 
\[
hA = \left\{ a_1+\cdots + a_h: a_i \in A \text{ for all } i \in [1,h] \right\}.
\]
It is a fundamental unsolved problem in additive number theory, and a problem 
that strangely seems not to have been previously investigated, 
to compute and understand the sumset sizes of finite sets 
of integers.\footnote{This 
 problem was discussed at three recent number theory conferences 
  (Integers 2025 in Athens, Georgia, CANT 2025 in New York, and BIRS-IMAG  in Granada, Spain) and no one could cite earlier work.  
  References to previous papers would be appreciated.}
For all positive integers $h$ and $k$, we define 
the \emph{range of sumset sizes}  
\[
\mcr_{\Z}(h,k) = \left\{ |hA|: A \subseteq \Z \text{ and } |A| = k \right\}.
\] 
This set was introduced in Nathanson~\cite{nath25bb} 
and studied in~\cite{nath25aa,nath2505,nath2505b}.

The idea of ``sumset'' is central in additive number theory 
and suggests many simple but \emph{a priori} difficult questions.  
For example, does there exist a set $A$ of size $|A| = 11$ 
whose $7$-fold sum has size $|7A| = 15551$?

The answer is ``yes,'' and the set 
\[
A_0 = \{ 0,897, 2056,2441,2988,3259,5294,6506,8013,9391,9872\}
\] 
has the desired sumset sum.  However, the answer was obtained by cheating: I started with the set $A_0$ and  computed its $7$-fold sum.  
For this set $A_0$, Maple quickly calculated the sumsets $hA_0$ 
and the  sumset sizes $|hA_0|$ 
for $h \in [1,7]$.  The second line of the following table 
gives the sequence of sumset sizes $|hA_0|$.  
The maximum size of an $h$-fold sumset of a set of size $k$ 
is $\binom{h+k-1}{h}$. 
The third line of the table gives the maximum size of an $h$-fold sumset 
of a set $A$ of size 11 for $h \in [1,7]$. 

\begin{center}
\begin{tabular}{c|r|r|r|r|r|r|r}
$h$ & 1 & 2 & 3 &4 & 5 & 6 & 7 \\ \hline
$|hA_0|$& 11 & 66 & 285 &  986 & 2878 & 7226 & 15551 \\ \hline 
$\max |hA|$ &  11 & 66 & 286 & 1001 & 3003 & 8008 & 19448 
\end{tabular}
\end{center}

 A $B_h$-set is a set $A$ of integers such that every integer 
 in the $h$-fold sumset $hA$ has a unique representation as a sum of 
 $h$ increasing elements of $A$.  
 If $A$ has size $k$, then $A$ is a $B_h$-set if and only if 
 $|hA| =  \binom{h+k-1}{h}$. 
 A $B_2$-set is also called a Sidon set.
 We see from the table above that $A_0$ is a Sidon set.  
Every $B_3$-set of size 11 satisfies $|3A| = 286$, but the special set $A_0$ 
has $|3A_0| = 285$, and so is not a $B_3$-set.  
Moreover, there must be exactly one integer $n_0$ that has two 
representations as a sum of three elements of $A_0$ and 
no integer that has more than two representations 
as a sum of three elements of $A_0$.  We have 
\[
n_0 = 18782 = 0 + 9391 + 9391 = 897 + 8013 + 9872 \in 3A_0.
\]

There are many related questions.  
The maximum size of a 7-fold sumset of a set of 10 integers is 
$\binom{16}{7} = 11440$.  It follows that if $A$ is a set of integers 
with $|7A| = 15551$, then $|A| \geq 11$.  
What is the cardinality of the  largest set  $A$ such that $|7A| = 15551$?

Sets $A$ and $B$ are \emph{affinely equivalent} if there are numbers 
$\lambda \neq 0$ and $\mu$ such that 
\[
B = \lambda\ast A + \mu = \{\lambda a+ \mu: a \in A\}.
\]
Affinely equivalent sets have the same sumset sums.  
Sets $A$ and $B$ are \emph{affinely inequivalent} 
if they are not affinely equivalent.  
How many affinely inequivalent sets $A$ satisfy $|A| = 11$ 
and $|7A| = 15551$?
How many affinely inequivalent finite sets $A$ satisfy $|7A| = 15551$?
The computation complexity inequality~\eqref{complexity} 
implies that all these questions can be answered in finite time.    

We have 
\[
\max \mcr_{\Z}(h,k) = \binom{h+k-1}{k-1} 
\]
and sets with this sumset size ($B_h$-sets and Sidon sets) have been extensively studied (\cite{nath-ubiq,nath-Canada,nath2502,obry04}).
Similarly, 
\[
\min \mcr_{\Z}(h,k) = h(k-1)+1
\]
and sets with this sumset size (called arithmetic progressions) 
and sets with very small sumset size have also been extensively studied 
(Freiman's theorem and related work~\cite{nath1996bb,tao-vu06}). 
Thus,
\beq      \label{R-interval}
\mcr_{\Z}(h,k) \subseteq \left[ h(k-1)+1, \binom{h+k-1}{k-1} \right].
\eeq

Because the set $\mcr_{\Z}(h,k)$ is finite, we can define the 
integer $N(h,k)$ as the smallest integer $N$ such that 
\[  
\mcr_{\Z}(h,k) = \{ |hA|: A \subseteq [0,N-1] \text{ and } |A| = k \}.
\]
The integer $N(h,k)$ determines the computational complexity of the 
sumset size problem.  It is proved in~\cite{nath2505b} 
that, for all $h \geq 3$ and $k \geq 3$,
\beq         \label{complexity}
N(h,k) < 4(8h)^{k-1}.
\eeq
Thus, for fixed $k$, the set $\mcr_{\Z}(h,k)$ can be computed in polynomial time.

Here are some simple results. 
For all $h$ and $k$, 
\[
\mcr_{\Z}(h,1) = \{1\}, \quad \mcr_{\Z}(1,k) = \{k\}, 
\qand
\mcr_{\Z}(h,2) = \{h+1\}.
\]
For all $k$, the set $\mcr_{\Z}(2,k)$ is the interval of integers:
\[
\mcr_{\Z}(2,k) = \left[ 2k-1, \binom{k+1}{2} \right]. 
\]
However, the set $\mcr_{\Z}(3,3)$ is not an interval of integers. We have 
\[
\mcr_{\Z}(3,3) \subseteq \left[ 7, 10 \right] 
\]
but 
\[
\mcr_{\Z}(3,3) = \{7,9,10\}.
\]
The number 8 missing.  For $k=3$ and all $h \geq 1$, 
Nathanson~\cite{nath25bb} proved that  
\[
\mcr_{\Z}(h,3) = \left\{  \binom{h+2}{2}  -  \binom{t}{2} : t \in [1,h] \right\}. 
\]
and so $|\mcr_{\Z}(h,3)| = h$.

\bprob
For  $k=4$, we have 
\[
\mcr_{\Z}(h,4) \subseteq \left[ 3h+1, \binom{h+3}{3} \right].
\]
Compute $\mcr_{\Z}(h,4)$ for all positive integers $h$.   
\eprob

\bprob
For all positive integers $h$ and $k$,  
compute the set of sumset sizes $\mcr_{\Z}(h,k)$.  
Explain the  integers ``missing'' from the interval 
on the right side of~\eqref{R-interval}.  
This is a core problem in additive number theory.  
\eprob

\bprob
 Let $k \geq 4$.  
Prove that ``most'' numbers in the interval 
on the right side of~\eqref{R-interval} 
are \emph{not} sumset sizes, that is, 
prove that 
\[
\left| \mcr_{\Z}(h,k) \right|= o_k \left( h^{k-1} \right)
\]
or, better, that 
\[
\left| \mcr_{\Z}(h,k) \right|= O_k\left( h^{k-2} \right). 
\]
\eprob

\bprob
Is there some kind of relation or reciprocity 
between the sets $\mcr_{\Z}(h,k)$ 
and $\mcr_{\Z}(k,h)$?  For example, is there a function $f: \N \rightarrow \N$ 
such that, if $t \in \mcr_{\Z}(h,k)$, then $f(t) \in \mcr_{\Z}(k,h)$? 
Is there a function $g: \N \rightarrow \N$ 
such that $t \in \mcr_{\Z}(h,k)$ if and only if $g(t) \in \mcr_{\Z}(k,h)$? 
Are there  functions $\alpha,\beta,\gamma$ such that, 
if $t \in \mcr_{\Z}(h,k)$, then $\gamma(t) \in \mcr_{\Z}(\alpha(h), \beta(k))$? 
\eprob

For related work on ``sumset size races,'' 
see~\cite{fox-krav-zhan25,krav25,peri-roto25}.

\section{5-fold sums of sets of size 4}
In this section we describe some computational experiments on 
5-fold sumsets of sets of size 4.  The calculations were programmed in Maple. 

\textbf{Experiment 1.} 
Let $q = 100$. 
There are 
\[
\binom{100}{4} = 3,921,225
\]
subsets of size 4 contained in the integer interval $[0,q-1]$.  
Their 5-fold sumsets are contained in $[16,56]$.  
For all $t \in [16,56]$, let $S(t)$ denote the number of sets 
$A \subseteq [0,q-1] = [0,99]$ such that  $|A| = 4$ and $|5A| = t$.  
The following table records the 
distribution of  sumset sizes.

\begin{center}
\begin{tabular}{cr|cr|cr} 
t & $S(t)$  & t & $S(t)$   & t& $S(t)$   \\ \hline
16& 1617 & 30&  4655 & 44 & 10574   \\
 17 & 0 & 31 & 3444& 45 & 18862 \\
 18 & 0 & 32 & 5273 & 46 & 450320 \\
 19 & 0 & 33 & 1906 &  47 &6706 \\
 20 & 2400  &  34& 5846 & 48& 6120  \\                      
21 & 1200  & 35 & 5400  &   49 &  0 \\
22 &0  & 36  & 280850  & 50 & 15002\\
23 & 0 & 37 & 1152  &  51 & 27332 \\
24& 4750  & 38& 3176  & 52 &  460734  \\
25 &  0 & 39  & 4612  &  53 & 7628 \\
26 &950 & 40  &  0 &  54 & 30688 \\
27 &4704 & 41  & 10012  &  55 & 667836 \\
28 & 0 & 42  & 5888  &  56 & 1867410 \\
29&  1568 &  43 &  2610 & &      \\
\end{tabular}
\end{center}
The relation $S(16) = 1617$ is equivalent to the statement that the interval 
$[0,99]$ contains exactly 1617  arithmetic progressions of length 4. 
The relation $S(56) = 1867410$ means that the interval 
$[0,99]$ contains exactly 1867410 $B_5$-sets.  

\textbf{Experiment 2.}
Let $q = 1,000$. 
Maple randomly chose 150,000 sets of size 4 
in the interval $[0,q-1]$ and computed their 5-fold sumsets.  

 The set of  experimentally observed sumset sizes is 
\begin{align*}
\{20, 27, 30,  32, 34, 35,  36,    39, 41, 42, 
  44, 45, 46, 47,   50, 51, 52, 53, 54, 55, 56\}. 
\end{align*}
The following table lists these elements $t$ of the set $R(5,4)$ and the number $S(t)$ of sets $A$  in the sample space with $|5A|= t$: \\

\begin{center}
\begin{tabular}{cr|cr}
$t$ & $S(t)$  & $|5A|$ & $S(t)$   \\ \hline
20 & 1 & 45 & 6 \\  
27 & 2 & 46 & 1738   \\                      
30 & 2  & 47& 1  \\
32& 2 &   50 & 6 \\
34& 2 &     51& 7   \\
35 & 1 &    52& 1957   \\
36& 982 &   53 & 2     \\
39& 2  & 54& 9  \\
41& 2 &  55 & 3251 \\
42 & 1 & 56& 142020 \\
44& 6 &    &  
\end{tabular}
\end{center}

 \textbf{Experiment 3.}
Let $q = 10,000$. 
Maple randomly chose 1,000,000 sets of size 4 
in the interval  $[0,q-1]$  and computed their 5-fold sumsets.  
 The set of  experimentally observed sumset sizes is  
\begin{align*}
\{21,36, 38, 39, 42, 44, 45, 46, 50, 51, 52, 54,  55, 56\}. 
\end{align*}
Let $q = 10,000$. 
The following table lists these elements $t$ of the set $R(5,4)$ 
and the number $S(t)$ of sets $A$  in the sample space with $|5A|= t$: \\

\begin{center}
\begin{tabular}{cr|cr}
$t$ & $S(t)$  & $t$ & $S(t)$ \\ \hline
21 & 584 & 46& 1429  \\
36 & 805 & 50 & 1 \\
38 & 1 &  51 & 1 \\
39 & 1 & 52 & 1654 \\
42 & 1 & 54 & 1 \\
44 &  1&  55 & 2559 \\
45 & 1 &  56 & 992961  
\end{tabular}
\end{center}

\textbf{Experiment 4.}
Let $q = 1,000,000$. 
Maple randomly chose 250,000 sets of size 4 
in the interval  $[0,q-1]$  and computed their 5-fold sumsets.  
 The set of  experimentally observed sumset sizes is  
\begin{align*}
\{ 36,  46,  52,   55, 56\}. 
\end{align*}
The following table lists these elements $t$ of the set $R(5,4)$ and the number $S(t)$ of sets $A$  in the sample space with $|5A|= t$: \\

\begin{center}
\begin{tabular}{cr }
$t$ & $S(t)$    \\ \hline
36 & 1 \\
 46& 4 \\
52 & 2\\
 55 &4 \\
 56 & 249989  
\end{tabular}
\end{center}

\bprob
In all four experiments, the sets with sumset size 56 predominate 
and these are the $B_5$-sets.  
The most popular sumset sizes are 56, 55, 52, 46, and 36. 
Their differences are the consecutive triangular numbers 1,3, 6, and 10.  
Why? 
\eprob

\bprob
Let $N = N(5,4) = 4\cdot 40^3 = 256,000$. Inequality~\eqref{complexity} 
implies that $\mcr_{\Z}(5,4)$ is the set of all 
5-fold sums of  4-element sets of integers contained 
in the interval $[0,N-1]$. 
Compute $\mcr_{\Z}(5,4)$. 

\eprob

\section{A mystery about popular sumset sizes}
Maple was used to generate random sets of 4 integers contained in the intervals $[0,999]$ and $[0,9999]$, to calculate 
their $h$-fold sumsets for $h \in [2,9]$, and to find the ``popular'' sumset sizes.  We obtained the following results.    \\

\begin{center}
\begin{tabular}{c|r|r}
$h$  & most popular sumset sizes & successive differences  \\ \hline
2 &      $   9, 10  $ & 1 \\ 
3 &      $    16,19, 20 $ &  3,1\\            
4 &     $    25,  31, 34, 35  $ & 6,3,1 \\
5 & $ 36, 46, 52, 55, 56 $ & 10, 6,3,1  \\
6 &  $49, 64, 74, 80, 83, 84 $& 15, 10, 6,3,1 \\
7 &   $ 64, 85, 100, 110, 116, 119, 120 $ & 21, 15, 10, 6,3,1   \\
8 & $   81, 109, 130, 145, 155, 161, 164, 165  $   & 28, 21, 15, 10, 6,3,1     \\
9 &     $   100, 136, 164, 185, 200, 210, 216, 219, 220 $   & 36, 28, 21, 15, 10, 6,3,1 \\  
\end{tabular}
\end{center}

\bprob
For all $h \in [2,9]$, the sumset size set $\mcr(h,4)$
 has $h$ popular sizes and the differences between consecutive 
 popular sumset sizes are the first $h-1$ 
triangular numbers.  Why? 
\eprob

Recall that the \emph{$j$th tetrahedral number} 
$\Tet_j = \binom{j+2}{3}$ is the sum of the first $j$ triangular numbers.  
The popular sumset sizes appear to be 
the differences between  $\binom{h+3}{3}$, 
which is the size of a 4-element 
$B_h$-set, and the $h$ consecutive tetrahedral numbers 
$\Tet_0,  \Tet_1,\ldots,  \Tet_{h-1}$. 

\bprob
Let $h \geq 3$.  
Prove that 
\[
\left\{ \binom{h+3}{3} - \Tet_j : j \in [1,h-1]\right\}  
\subseteq \mcr_{\Z}(h,4). 
\]
\eprob

\bprob
Let $h$ and $q$ be positive integers.  
Let 
\[
\mca(4,q) = \left\{ A \subseteq [0,q-1]: |A| = 4  \right\} 
\]
\[
\mca^*(4,q,h) = \left\{ A \in \mca(4,q):
\text{$A$ is not a $B_h$-set} \right\} 
\]
and 
\[
\mcb^*(4,q,h ) = \left\{ A \in \mca^*(4,q,h ) :
 |hA| \in \left\{ \binom{h+3}{3} - \Tet_j : j \in [1,h-1]\right\}  \right\} 
\]
Almost all 4-element subsets of $[0,q-1]$ are $B_h$-sets (Nathanson~\cite{nath-ubiq,nath-Canada}):
\[
\lim_{q\rightarrow \infty} \frac{| \mca^*(4,q,h) |}{|\mca(4,q) |} = 0.
\]
Is it true that the sumsets sizes of almost all sets $A \in \mca(4,q)$ that 
are not $B_h$-sets  (that is, sets in $\mca^*(4,q)$) are the popular sumset sizes?  Equivalently, is it true that 
\[
\lim_{q\rightarrow \infty} \frac{| \mcb^*(4,q,h) |}{|\mca^*(4,q,h) |} = 1.
\]
\eprob

\bprob
Discover the pattern of popular sumset sizes for 
sets of size $k \geq 5$. 
\eprob

\section{A toy problem for physics} 
For positive integers $h$ and $k$, let $\mcx_{h,k}$ be the 
set of $k$-tuples  $\mbx = (x_1,\ldots, x_k)$ of nonnegative integers 
such that $\sum_{i=1}^k x_i = h$.  
Let $A = \{a_1,\ldots, a_k\}$ be a set of $k$ integers 
and let $\vec{A} = (a_1,\ldots, a_k)$ be the associated lattice point 
in $\Z^n$.  
For all $\mbx = (x_1,\ldots, x_k) \in \mcx_{h,k}$, we define 
\[
\mbx\cdot \vec{A} = \sum_{i=1}^k x_ia_i. 
\]
Then 
\beq                \label{hAX}
hA = \left\{\mbx\cdot \vec{A}: \mbx \in \mcx_{h,k} \right\}. 
\eeq
Because $\mcx_{h,k}$ is invariant with respect to the permutation 
group $S_n$ (that is, $\mbx = (x_1,\ldots, x_k) \in \mcx_{h,k}$  
if and only if $\sigma \mbx = (x_{\sigma(1)},\ldots, x_{\sigma(k)}) \in \mcx_{h,k}$ 
for all $\sigma \in S_n$),  relation~\eqref{hAX} 
is independent of the ordering of the coordinates in the vector $\vec{A}$.

The \emph{representation function} $r_{A,h}(n)$ counts the number 
of representations of $n$ as a sum of $h$ elements of $A$, that is, 
the number of $k$-tuples $\mbx = (x_1,\ldots, x_k) \in \mcx_{h,k}$ 
such that $\mbx \cdot \vec{A} = \sum_{i=1}^k x_ia_i = n$.  
Analogous to collisions of molecules in a gas, 
$k$-tuples $\mbx  \in \mcx_{h,k}$ and $\mby  \in \mcx_{h,k}$ 
``collide'' with respect to the set $A$ if the $h$-fold sums 
$\mbx\cdot \vec{A} $ and $\mby \cdot \vec{A}$ are equal.  
Moreover, like multiple collisions of molecules, 
multiple collisions of vectors in $\mcx_{h,k}$, that is, integers $n \in hA$ 
such that $r_{A,h}(n) \geq 3$, should  usually 
be much rarer than simple collisions with $r_{A,h}(n) = 2$, 
but not always.   For example, with $k \geq 5$, 
the integer interval $A = [0,k-1]$ 
satisfies $2A = [0, 2k-2]$ and  
\[
r_{A,2}(n) = \begin{cases}
\lfloor \frac{n}{2}\rfloor + 1& \text{if $0 \leq n \leq k-1$}\\
\lfloor \frac{2k-2-n}{2}\rfloor + 1& \text{if $k-1 \leq n \leq 2k-2$.} 
\end{cases}
\]
Thus,
\begin{align*}
r_{A,2}(n) & = 1      \qquad  \text{if $n \in \{ 0,1,2k-3,2k-2 \}$} \\
r_{A,2}(n) & =  2     \qquad  \text{if   $n \in \{ 2,3 ,2k-5,2k-4 \}$} \\  
r_{A,2}(n) & \geq 3 \qquad   \text{if  $4 \leq n \leq 2k-6$  } 
\end{align*}
and so 
\beq       \label{r3>r2}
\card\left( n \in \Z: r_{A,h}(n) \geq 3 \right) 
> \card\left( n \in \Z: r_{A,h}(n) = 2 \right)
\eeq
for all $k \geq 7$. 
Inequality~\eqref{r3>r2} should be atypical.  
A quantitative version of this remark is the following.

\bprob
For every sufficiently large integer $q$, let 
\[
\mca(k,q )= \left\{ A \subseteq [0,q-1]: |A| = k\right\}.
\]
Prove that 
\[
\card\left( n \in \Z: r_{A,h}(n) \geq 3 \right) 
< \card\left( n \in \Z: r_{A,h}(n) = 2 \right) 
\]
for almost all $A \in \mca(k,q)$. 
\eprob

Computing the patterns of sumset sizes of finite sets of integers can be considered as a toy problem for certain much deeper questions in 
mathematical physics, such as the Hilbert program of deducing 
macroscopic fluid mechanical phenomena from 
mesoscopic  and microscopic motions.  
A significant part  of recent work of Yu Deng,  
Zaher Hani, and Xiao Ma~\cite{deng-hani-ma25,slom25} 
on Hilbert's sixth problem (the derivation of the Navier-Stokes equation 
from the Boltzman  equation from elastic collisions of a gas) 
is combinatorial and has a number theoretical flavor.  
It would be of interest to examine the combinatorial results of 
Deng, Hani, and Ma, and and the similarities between phenomena 
of combinatorial additive number theory and fluid dynamics.

\emph{Acknowledgements.}
I thank Christian Elsholtz, Harald Helfgott, Noah Kravitz, Kevin O'Bryant, 
and Steven Senger for helpful remarks on this paper.

\end{document}